\documentclass[12pt,oneside,reqno]{amsart}
\usepackage[margin=1in]{geometry}
\usepackage[utf8]{inputenc}
\usepackage{amsthm, amssymb}
\usepackage[colorlinks=true, pdfstartview=FitV, linkcolor=blue, citecolor=blue, urlcolor=blue]{hyperref}
\usepackage[colorinlistoftodos]{todonotes}
\usepackage{algorithm}
\usepackage{algpseudocode}
\usepackage{pgf}
\usepackage{tikz}
\usepackage{tikz-cd}
\usetikzlibrary{positioning,shapes,shadows,arrows}
\usepackage{bbm,bm}
\usepackage{enumitem}
\usepackage{mathabx}
\usepackage{comment}
\usepackage{ytableau}
\usepackage{thmtools}
\usepackage{thm-restate}

\newtheorem{prop}{Proposition}[section]
\newtheorem{thm}[prop]{Theorem}

\newtheorem{lemma}[prop]{Lemma}
\newtheorem{cor}[prop]{Corollary}
\newtheorem{conj}[prop]{Conjecture}
\theoremstyle{remark}
\newtheorem{exa}[prop]{Example}
\newtheorem{rem}[prop]{Remark}
\newtheorem{defn}[prop]{Definition}

\newcommand{\topGro}{\widehat{\fG}}

\newcommand{\fG}{\mathfrak{G}}

\newcommand{\wt}{\mathsf{wt}}

\newcommand{\PD}{\mathsf{PD}}

\newcommand{\BVPD}{\mathsf{BVPD}}
\newcommand{\MVPD}{\mathsf{MVPD}}

\newcommand{\x}{\textbf{x}}
\newcommand{\y}{\textbf{y}}
\newcommand{\wty}{\mathsf{wty}}
\newcommand{\raj}{\mathsf{raj}}
\newcommand{\maj}{\mathsf{maj}}
\newcommand{\droop}{\mathsf{droop}}
\newcommand{\Supp}{\mathsf{Supp}}

\newcommand{\btile}{
 \begin{tikzpicture}[x=0.8em,y=0.8em,thick,color = blue]
\draw[step=1,gray,thin] (0,0) grid (1,1);
\draw[color=black, thick, sharp corners] (0,0) rectangle (1,1);
\end{tikzpicture}}

\newcommand{\htile}{
\begin{tikzpicture}[x=0.8em,y=0.8em,thick,color = blue]
\draw[step=1,gray,thin] (0,0) grid (1,1);
\draw[color=black, thick, sharp corners] (0,0) rectangle (1,1);
\draw(1.0,0.5)--(0.0,0.5);
\end{tikzpicture}
}

\newcommand{\vtile}{
\begin{tikzpicture}[x=0.8em,y=0.8em,thick,color = blue]
\draw[step=1,gray,thin] (0,0) grid (1,1);
\draw[color=black, thick, sharp corners] (0,0) rectangle (1,1);
\draw(0.5,1.0)--(0.5,0.0);
\end{tikzpicture}}

\newcommand{\ptile}{
\begin{tikzpicture}[x=0.8em,y=0.8em,thick,color = blue]
\draw[step=1,gray,thin] (0,0) grid (1,1);
\draw[color=black, thick, sharp corners] (0,0) rectangle (1,1);
\draw(0.5,1.0)--(0.5,0.0);
\draw(1.0,0.5)--(0.0,0.5);
\end{tikzpicture}}

\newcommand{\rtile}{
\begin{tikzpicture}[x=0.8em,y=0.8em,thick,rounded corners,color = blue]
\draw[step=1,gray,thin] (0,0) grid (1,1);
\draw[color=black, thick, sharp corners] (0,0) rectangle (1,1);
\draw(1.0,0.5)--(0.5,0.5)--(0.5,0.0);
\end{tikzpicture}}

\newcommand{\jtile}{
\begin{tikzpicture}[x=0.8em,y=0.8em,thick,rounded corners,color = blue]
\draw[step=1,gray,thin] (0,0) grid (1,1);
\draw[color=black, thick, sharp corners] (0,0) rectangle (1,1);
\draw(0.5,1.0)--(0.5,0.5)--(0.0,0.5);
\end{tikzpicture}}

\newcommand{\mtile}{
\begin{tikzpicture}[x=0.8em,y=0.8em,thick,rounded corners,color = blue]
\draw[step=1,gray,thin] (0,0) grid (1,1);
\draw[color=black, thick, sharp corners] (0,0) rectangle (1,1);
\draw(1.0,0.5)--(0.5,0.5)--(0.5,0.0);
\node at (0.5,0.5) {$\bullet$};
\end{tikzpicture}}

\newcommand{\bumptile}{
\begin{tikzpicture}[x=0.8em,y=0.8em,thick,rounded corners,color = blue]
\draw[step=1,gray,thin] (0,0) grid (1,1);
\draw[color=black, thick, sharp corners] (0,0) rectangle (1,1);
\draw(1.0,0.5)--(0.5,0.5)--(0.5,0.0);
\draw(0.5,1.0)--(0.5,0.5)--(0.0,0.5);
\end{tikzpicture}}

\newcommand{\bumphtile}{
\begin{tikzpicture}[x=0.8em,y=0.8em,thick,rounded corners,color = blue]
\draw[step=1,gray,thin] (0,0) grid (2,1);
\draw[color=black, thick, sharp corners] (0,0) rectangle (2,1);
\draw(2.0,0.5)--(0.5,0.5)--(0.5,0.0);
\draw(0.5,1.0)--(0.5,0.5)--(0.0,0.5);
\end{tikzpicture}}

\newcommand{\rhtile}{
\begin{tikzpicture}[x=0.8em,y=0.8em,thick,rounded corners,color = blue]
\draw[step=1,gray,thin] (0,0) grid (2,1);
\draw[color=black, thick, sharp corners] (0,0) rectangle (2,1);
\draw(2.0,0.5)--(0.5,0.5)--(0.5,0.0);
\end{tikzpicture}}

\definecolor{darkblue}{rgb}{0.0,0,0.7} 
\definecolor{darkred}{rgb}{0.7,0,0} 
\definecolor{darkgreen}{rgb}{0, .6, 0} 
\newcommand{\definition}[1]{{\color{darkred}\emph{#1}}} 

\title{Grothendieck polynomials of inverse fireworks permutations}

\author[C.~Chou]{Chen-An (Jack) Chou}
\address[C. Chou]{Department of Mathematics, University of Florida, Gainesville, Fl 32611, U.S.A.}
\email{c.chou@ufl.edu}

\author[T.~Yu]{Tianyi Yu}
\address[T. Yu]{Laboratoire d'Alg\`ebre, de Combinatoire et d'Informatique Math\'ematique, Universit\'e du Qu\'ebec \`a Montr\'eal, Montreal QC, Canada}
\email{yu.tianyi@uqam.ca}

\begin{document}
\maketitle
\begin{abstract}
Pipe dreams are combinatorial objects that compute Grothendieck polynomials.
We introduce a new combinatorial object that naturally 
recasts the pipe dream formula.
From this, we obtain the first direct 
combinatorial formula for 
the top degree components of Grothendieck polynomials,
also known as the Castelnuovo-Mumford polynomials.
We also prove the inverse fireworks case of a conjecture of 
M{\'e}sz{\'a}ros, Setiabrata, and St. Dizier
on the support of Grothendieck polynomials. 
\end{abstract}

\section{Introduction}
\label{S: Intro}
Fix $n \in \mathbb{Z}_{>0}$ throughout the paper. 
For $w \in S_n$,
Lascoux and Sch\"utzenberger \cite{LS:Groth} introduced 
the Grothendieck polynomials $\fG_w(\x)$, which are explicit polynomial representatives of the 
$K$-classes of structure sheaves of Schubert varieties in flag varieties.
In general, the Grothendieck polynomials are not homogeneous. 
Their lowest degree homogeneous components are the Schubert polynomials, 
which represent the cohomology classes of Schubert varieties in flag varieties.

Grothendieck polynomials can be computed combinatorially using objects called \definition{pipe dreams (PD)}~\cite{BB, BJS, FK} (see Definition~\ref{D: PD}). 
Each $w \in S_n$
is associated with a set of PDs,
denoted as $\PD(w)$.
Each $P \in \PD(w)$
has a set of weight tiles $\wty(P)$,
which leads to a monomial $\wt_P(\x)$.
By Fomin and Kirillov~\cite{FK}, one may compute the Grothendieck polynomial of $w$ as
\begin{align*}
    \textrm{$\fG$}_w(\x) = \sum_{P \in \PD(w)}  (-1)^{|\wty(P)| - \ell(w)}\wt_P(\x).
\end{align*}
We treat this equation as our definition of $\fG_w(\x)$.
The main innovation of this paper
is an alternative perspective on PDs. 
We remove certain pipes from
a PD to obtain a novel
combinatorial object called a
\definition{marked vertical-less pipe dream (MVPD)} (see
Definition~\ref{D: MVPD}),
while retaining exactly the same information.
In other words, we introduce a set $\MVPD(w)$
and a bijection $\Phi: \PD(w) \rightarrow \MVPD(w)$.
We also define the weighty tiles $\wty(M)$
of each $M \in \MVPD(w)$ so that $\Phi$ preserves $\wty(\cdot)$.
Consequently, each $M \in \MVPD(w)$
is associated with a monomial $\wt_M(\x)$
which agrees with its corresponding PD. 
This set $\MVPD(w)$ 
recasts the PD formula in a way
that helps us derive the following
two applications on Grothendieck polynomials. 

\subsection{Combinatorial formula
for Castelnuovo-Mumford polynomials}
There has been a growing recognition of the significance of the Castelnuovo-Mumford polynomial $\topGro_w(\x)$, 
which is defined as
the top degree homogeneous component of $\fG_w(\x)$. 
For instance, 
\begin{itemize}
\item the degree of $\topGro_w(\x)$
determines the Castelnuovo–Mumford regularity of matrix Schubert varieties~\cite{RRRSW};
\item the support of $\topGro_w(\x)$ conjecturally governs
the support of $\fG_w(\x)$~\cite[Conjecture 1.3]{MSS}.
\end{itemize}

With these motivations, 
there has been a recent surge in the study of $\topGro_w(\x)$~\cite{CY, DMS, Haf, HMSS, PSW}.
Notably, Pechenik, Speyer, and Weigandt~\cite{PSW} showed
that each 
$\topGro_w(\x)$ is an integer multiple
of $\topGro_u(\x)$ for some inverse fireworks permutation $u$, defined as follows. 

\begin{defn}
Any permutation can be decomposed 
into decreasing runs.
A permutation is called \definition{fireworks} if the first
number in each decreasing run is increasing.
A permutation is called an \definition{inverse fireworks permutation} if its inverse is fireworks.    
For instance, $\textbf{3}1\textbf{67}542$ is fireworks while 
$\textbf{6}1\textbf{37}542$ is not.
\end{defn}
 
Thus, to understand all $\topGro_w(\x)$,
one might focus on $\topGro_u(\x)$
for inverse fireworks $u$.
As far as the authors are aware, 
there does not exist any direct combinatorial formula for $\topGro_w(\x)$,
besides extracting highest degree elements
from $\PD(w)$.
Ideally, for a polynomial $f_w$ indexed by
$w \in S_n$,
one would wish a combinatorial formula
following the template:
\begin{quote}
The polynomial
$f_w$ is a sum over fillings of an $n \times n$
grid
using tile set $\textbf{T}$, such that when we trace
the pipes, it satisfies some boundary condition associated with $w$.
Each such filling contributes a monomial,
whose power of $x_i$ is the number of tiles 
in $\textbf{T'} \subseteq \textbf{T}$ in row $i$.    
\end{quote}

The pipe dream for 
$\fG_w$ formula is one such formula, 
where $\textbf{T}$ consists of $\ptile$, $\bumptile$, $\jtile$ and $\btile$.
The marked bumpless pipe dream formula for $\fG_w$~\cite{LLS, LLS2}
is also one such formula. 
However, Schubert polynomials do not have a formula fitting into this template.
Its pipe dream formula requires an extra condition: 
Two pipes cannot cross more than once.

We provide such an ideal combinatorial formula 
of $\topGro_w(\x)$ for inverse fireworks $w$.
Our $\textbf{T}$ consists of five tiles:
$$
\btile,\quad \quad \htile, \quad \quad \ptile, \quad \quad \jtile, \quad \quad \rtile
$$
Since $\textbf{T}$ does not have 
$\bumptile$ and $\vtile$, 
we call our objects 
bumpless vertical-less pipe dreams (BVPDs).
We prove our formula by establishing a bijection 
between BVPDs and the highest degree elements from $\MVPD(w)$.
Via the bijection between $\MVPD(w)$ and $\PD(w)$, 
we also obtain a characterization of 
highest degree PDs of inverse fireworks permutations.  

\subsection{Conjectures on the supports of Grothendieck polynomials}

Huh, Matherne, M{\'e}sz{\'a}ros and 
St. Dizier~\cite{HMMS}
conjectured
that homogenized Grothendieck polynomials are Lorentzian 
(up to appropriate normalization).
Their conjectures inspired 
recent studies on the support of $\fG_w$, 
which is the set of monomials that appear
in $\fG_w$.
We denote it as $\Supp(\fG_w)$.
Notably, M{\'e}sz{\'a}ros, Setiabrata, and St. Dizier~\cite{MSS} made the following conjectures.

\begin{conj}{\cite[Conj. 1.1]{MSS}}
\label{C: one}
    Let $m \in \Supp(\fG_w(\x))$. If the degree of $m$ is less than the degree of $\fG_w(\x)$, then there exists a distinct $m' \in \Supp(\fG_w(\x))$ such that $m$ divides $m'$.
\end{conj}

\begin{conj}{\cite[Conj. 1.2]{MSS}}
\label{C: main}
Let $m \in \Supp(\fG_w(\x))$. If the degree of $m$ is less than the degree of $\fG_w(\x)$, then there exists $i$ such that $m x_i \in \Supp(\fG_w(\x))$.
\end{conj}

Conjecture~\ref{C: main} is a natural strengthening of Conjecture~\ref{C: one}. 
M{\'e}sz{\'a}ros, Setiabrata, and St. Dizier~\cite{MSS} showed that Conjecture~\ref{C: one} holds for fireworks permutations. Hafner~\cite{Haf} proved Conjecture~\ref{C: main} for vexillary permutations. Recently, Chen, Fan, and Ye~\cite{CFY} proved Conjecture~\ref{C: main} for zero-one Grothendieck polynomials. We prove Conjecture~\ref{C: main} for inverse fireworks permutations. 

We may translate Conjecture~\ref{C: main} combinatorially:
For each $M \in \MVPD(w)$ 
that is not a top degree element,
we may find $M' \in \MVPD(w)$ such that 
$\wt_{M'}(\x) = \wt_M(\x) x_i$ for some $i$.
We prove Conjecture~\ref{C: main}
for an inverse fireworks permutation $w$ 
constructively: 
For each $M \in \MVPD(w)$,
we give an explicit algorithm that constructs $M' \in \MVPD(w)$.
The advantage of using 
$\MVPD(w)$ instead of $\PD(w)$,
is that we are allowed to perform the droop move,
which was originally defined on bumpless
pipe dreams~\cite{LLS}.
\begin{rem}
Anna Weigandt kindly informed us in private communication that she has an independent proof of Conjecture~\ref{C: main} for inverse fireworks permutations using a different technique.
\end{rem}

The rest of the paper is structured as follows.
In~\S\ref{S: Background},
we cover necessary background 
regarding pipe dreams and $\fG_w(\x)$.
In~\S\ref{S: MVPD},
we introduce marked vertical-less pipe dreams
and use them to rephrase the PD formula.
In~\S\ref{S: BVPD}, 
we introduce bumpless vertical-less 
pipe dreams and give a direct formula of $\topGro_w(\x)$.
In~\S\ref{S: Proving conjecture},
we prove Conjecture~\ref{C: main}
for inverse fireworks permutations by construction. 

\section{Background}
\label{S: Background}
\begin{defn}
\label{D: PD}
A \definition{pipe dream} (PD)~\cite{BB, BJS, FK}
is a tiling of an $n \times n$ grid,
using $\ptile, \bumptile, \jtile$ and $\btile$.
Pipes enter from the left and exit from the top.
We trace the pipes in a PD from left to top as follows. 
The pipe entering from row $p$ is called pipe $p$.
Suppose we see a $\ptile$ where the pipe
on the left (resp. bottom) has label $p$ (resp. $q$).
If pipe $p$ and $q$ have not crossed before, 
we say they cross in this tile
and let pipe $p$ (resp. $q$)
exit from the right (resp. top).
Otherwise, 
we let pipe $p$ (resp. $q$)
exit from the top (resp. right).
Notice that this rule is the same
as saying pipe $\max(p,q)$ exits
from the top and the other exits from the right.

After tracing the pipes, 
we may read off the labels of the pipes
on the top edge of the PD
as a permutation $w \in S_n$.
We say this PD is associated with $w^{-1}$.
Let $\PD(w)$ be the set of the pipe dreams
associated with $w$.
\end{defn}

\begin{exa}
The following is a pipe dream of the permutation $w$ with one-line notation $24513$. Its inverse has one-line notation $41523$.
$$
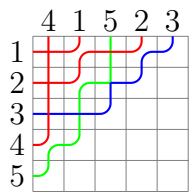
\begin{tikzpicture}[x=1em,y=1em,thick,rounded corners, color = blue]
\node[color=black] at (-0.5,4.5) {$1$};
\node[color=black] at (-0.5,3.5) {$2$};
\node[color=black] at (-0.5,2.5) {$3$};
\node[color=black] at (-0.5,1.5) {$4$};
\node[color=black] at (-0.5,0.5) {$5$};
\node[color=black] at (0.5,5.5) {$4$};
\node[color=black] at (1.5,5.5) {$1$};
\node[color=black] at (2.5,5.5) {$5$};
\node[color=black] at (3.5,5.5) {$2$};
\node[color=black] at (4.5,5.5) {$3$};
\draw[step=1,gray,ultra thin] (0,0) grid (5,5);
\draw[color=green] (0,0.5)--(0.5, 0.5)--(0.5, 1.5)--(1.5, 1.5)--(1.5, 3.5)--(2.5, 3.5);
\draw[color=green] (2.5, 3.5)--(2.5, 5);
\draw[color=blue] (2.5, 3.5)--(3.5, 3.5)--(3.5, 4.5)--(4.5, 4.5)--(4.5, 5);
\draw[color=blue] (0,2.5)--(2.5, 2.5)--(2.5, 3.5);
\draw[color=red] (0,1.5)--(0.5, 1.5)--(0.5, 5);
\draw[color=red] (0,3.5)--(1.5, 3.5)--(1.5, 4.5)--(3.5,4.5)--(3.5,5);
\draw[color=red] (0,4.5)--(1.5, 4.5)--(1.5, 5);
\end{tikzpicture}
$$
We make pipe $3$ blue and pipe $5$ green.
Notice that 
pipe $3$ and pipe $5$ cross at $(3,2)$.
However, pipe $3$ and pipe $5$ do not 
cross at $(2,3)$
since they already crossed. 
\end{exa}

For a PD $P$,
let $\wty(P)$ be the set of $(i,j)$ that is $\ptile$
in $P$.
Define 
$$
\wt_P(\x) = \prod_{(i,j) \in \wty(P)} x_i, \quad \wt_P(\x, \y) = \prod_{(i,j) \in \wty(P)} (x_i + y_j - x_i y_j).
$$

For $w \in S_n$,
let $\ell(w)$ be the number of $i < j$
such that $w(i) > w(j)$.
Following~\cite{FK} and~\cite{KM},
the \definition{Grothendieck polynomial} $\fG_w(\x)$
and the \definition{double Grothendieck polynomial} $\fG_w(\x, \y)$
can be defined as

\begin{align}
\label{EQ: Gro formulas}
\begin{split}
\fG_w(\x) & := \sum_{P \in \PD(w)}  (-1)^{|\wty(P)| - \ell(w)}\wt_P(\x),\\
\fG_w(\textbf{x}, \textbf{y}) & := \sum_{P \in \PD(w)} (-1)^{|\wty(P)| - \ell(w)} \wt_P(\x, \y).
\end{split}
\end{align}

\begin{exa}
The following are all the pipe dreams of the permutation $w=2413$:
$$
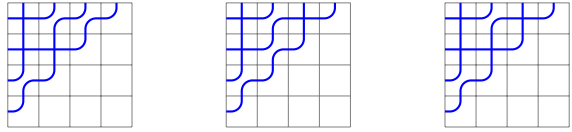
\begin{tikzpicture}[x=1em,y=1em,thick,rounded corners, color = blue]
\draw[step=1,gray,ultra thin] (0,0) grid (4,4);
\draw[color=blue] (0,0.5)--(0.5, 0.5)--(0.5, 1.5)--(1.5, 1.5)--(1.5, 3.5)--(2.5, 3.5)--(2.5, 4);
\draw[color=blue] (0,1.5)--(0.5, 1.5)--(0.5, 4);
\draw[color=blue] (0,2.5)--(2.5, 2.5)--(2.5, 3.5)--(3.5, 3.5)--(3.5, 4);
\draw[color=blue] (0,3.5)--(1.5, 3.5)--(1.5, 4);
\end{tikzpicture}
\quad\quad\quad
\begin{tikzpicture}[x=1em,y=1em,thick,rounded corners, color = blue]
\draw[step=1,gray,ultra thin] (0,0) grid (4,4);
\draw[color=blue] (0,0.5)--(0.5, 0.5)--(0.5, 1.5)--(1.5, 1.5)--(1.5, 2.5)--(2.5, 2.5)--(2.5, 4);
\draw[color=blue] (0,1.5)--(0.5, 1.5)--(0.5, 4);
\draw[color=blue] (0,2.5)--(1.5, 2.5)--(1.5, 3.5)--(3.5, 3.5)--(3.5, 4);
\draw[color=blue] (0,3.5)--(1.5, 3.5)--(1.5, 4);
\end{tikzpicture}
\quad\quad\quad
\begin{tikzpicture}[x=1em,y=1em,thick,rounded corners, color = blue]
\draw[step=1,gray,ultra thin] (0,0) grid (4,4);
\draw[color=blue] (0,0.5)--(0.5, 0.5)--(0.5, 1.5)--(1.5, 1.5)--(1.5, 3.5)--(3.5, 3.5)--(3.5, 4);
\draw[color=blue] (0,1.5)--(0.5, 1.5)--(0.5, 4);
\draw[color=blue] (0,2.5)--(2.5, 2.5)--(2.5, 3.5)--(2.5, 4);
\draw[color=blue] (0,3.5)--(1.5, 3.5)--(1.5, 4);
\end{tikzpicture}
$$
\end{exa}
Therefore,
\begin{align*}
    \fG_w(\x) &= x_1x_2^2+x_1^2x_2-x_1^2x_2^2 \\
    \fG_w(\x, \y) &= (x_1+y_1-x_1y_1)(x_2+y_1-x_2y_1)(x_2+y_2-x_2y_2)\\
    &+(x_1+y_1-x_1y_1)(x_2+y_1-x_2y_1)(x_1+y_3-x_1y_3)\\
    &-(x_1+y_1-x_1y_1)(x_2+y_1-x_2y_1)(x_2+y_2-x_2y_2)(x_1+y_3-x_1y_3)
\end{align*}

We make one simple observation on PDs
that will be useful later.
\begin{lemma}
\label{L: abc}
Say three pipes
enter a row of a PD from the bottom:
Pipe $a$ enters on the left of pipe $b$
and pipe $b$ enters on the left of pipe $c$.
Suppose pipe $a$ and pipe $b$
have not crossed, but pipe $a$
and pipe $c$ have crossed. 
Then pipe $b$ and pipe $c$
must have crossed.
\end{lemma}
\begin{proof}
Since pipe $a$ enters the row on
the left of pipe $b$ and they did not cross,
we have $a < b$.
Since pipe $a$ enters the row on
the left of pipe $c$ and they have crossed,
we have $c < a$.
Thus, $c < a < b$.
Since pipe $b$ enters the row on
the left of pipe $c$, they have crossed.
\end{proof}

The degree of $\fG_w(\x)$ was given by the 
statistic $\raj(w)$ defined by 
Pechenik, Speyer, and Weigandt~\cite{PSW}.
In other words, among PDs in $\PD(w)$,
the maximal number of weighty tiles is $\raj(w)$.
We define $\topGro_w(\x)$ as
$\sum_{P} \wt_P(\x)$
where the sum is over all $P \in \PD(w)$
with $\raj(w)$ weighty tiles. 
Up to a sign, $\topGro_w(\x)$
agrees with the 
top degree component of $\fG_w(\x)$.
In this paper, 
we need the following two properties 
of the $\raj$ statistic. 
\begin{prop}[{\cite[Proposition 3.8]{PSW}}]
\label{P: raj = maj}
For $w \in S_n$,
$\raj(w) = \maj(w)$ if and only if 
$w$ is fireworks. 
Here, 
$\maj(w)$ is the \definition{major index}
of $w$, defined as $\sum_{i: w(i) > w(i+1)} i$.
\end{prop}

\begin{cor}[{\cite[Corollary 4.5]{PSW}}]
\label{C: raj = raj inv}
For $w \in S_n$, $\raj(w) = \raj(w^{-1})$.
\end{cor}

\section{Marked Vertical-less Pipe Dreams}
\label{S: MVPD}

We introduce
combinatorial objects which 
we call marked vertical-less pipe dreams (MVPD).
An MVPD can be obtained by removing certain pipes from a PD. 
We rephrase~(\ref{EQ: Gro formulas}) and obtain MVPD formulas for $\fG_w(\x)$ and $\fG_w(\x, \y)$ in Corollary~\ref{C: MVPD}.

\begin{defn}
\label{D: MVPD}
A \definition{vertical-less pipe dream} (VPD) consists of the following six tiles:
$$
\btile,\quad \quad \htile, \quad \quad \ptile, \quad \quad \jtile, \quad \quad \rtile, \quad \quad \bumptile,
$$  
on an $n \times n$ grid.
Notice that we are not 
using the vertical tile $\vtile$.
The pipes of a VPD enter from the left edge
of the $n \times n$ grid
and exit from the top edge.
We trace pipes from left to top
in the same way as PDs.
The pipe entering from row $p$ is called pipe $p$.
A \definition{marked vertical-less pipe dream} (MVPD)
is a VPD where some $\rtile$
are marked as $\mtile$.
The pipe in a marked tile must have a 
$\htile$ in a lower row than the marked tile.
\end{defn}

The \definition{column-to-row code} of a MVPD
$M$ is a sequence of $n$ numbers. 
If there is no pipe exiting at column $c$ of $M$,
then the $c\textsuperscript{th}$ entry is $0$.
Otherwise, 
say pipe $r$ exits in column $c$,
then the $c\textsuperscript{th}$ entry is $r$.
When drawing a MVPD, we omit blank rows
on the bottom and blank
columns to the right. 

\begin{exa}
\label{E: MVPD}
Suppose $n = 8$.
The following is a MVPD which has column-to-row code \\
$(0,0,0,4,0,3,0,6)$. Notice that $(1,2),(2,1),$ and $(5,1)$ cannot be marked while $(3,4)$ may or may not be marked.
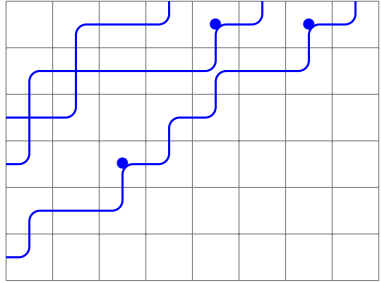
\begin{center}
    \begin{tikzpicture}[x=1.5em,y=1.5em,thick,rounded corners,color = blue]
    \draw[step=1,gray,ultra thin] (0,0) grid (8,6);
     \draw[color=blue, thick] (0,3.5)--(1.5,3.5)--(1.5,5.5)--(3.5,5.5)--(3.5,6);
     \draw[color=blue, thick] (0,2.5)--(0.5,2.5)--(0.5,4.5)--(4.5,4.5)--(4.5,5.5)--(5.5,5.5)--(5.5,6);
     \draw[color=blue, thick] (0,0.5)--(0.5,0.5)--(0.5,1.5)--(2.5,1.5)--(2.5,2.5)--(3.5,2.5)--(3.5,3.5)--(4.5,3.5)--(4.5,4.5)--(6.5,4.5)--(6.5,5.5)--(7.5,5.5)--(7.5,6);
     \node at (2.5,2.5) {$\bullet$};
     \node at (4.5,5.5) {$\bullet$};
     \node at (6.5,5.5) {$\bullet$};
    \end{tikzpicture}
\end{center}
\end{exa}

For a MVPD $M$,
we let $\wty(M)$ be the set of $(i,j)$
that is $\htile$, $\ptile$ or $\mtile$ in $M$.
We define 
$$
\wt_M(\x) = \prod_{(i,j) \in \wty(M)} x_i, 
\quad \wt_M(\x, \y) = \prod_{(i,j) \in \wty(M)} (x_i + y_j - x_i y_j).
$$

We associate certain MVPDs
to each permutation $w \in S_n$.
The \definition{left-to-right maxima}
of $w \in S_n$ are the numbers $w(i)$
such that $w(j) < w(i)$ for all $j < i$.
For instance, the left-to-right maxima
of the permutation with one-line notation $2143$
are $2$ and $4$.

\begin{rem}
\label{R: red pipe no x}
Notice that in the one-line notation of a permutation,
its left-to-right maxima must increase from left to right.
Consequently, in $P \in \PD(w)$, pipes labeled by left-to-right maxima of $w^{-1}$
cannot cross. 
\end{rem}

Take $w \in S_n$.
We start from $(w^{-1}(1), \cdots, w^{-1}(n))$
and turn the left-to-right maxima
of $w^{-1}$ into $0$.
Let $\alpha'(w)$ be the resulting sequence.
Finally, define $\MVPD(w)$
as the set of MVPDs with column-to-row
code $\alpha'(w)$.

\begin{exa}
Take $w \in S_n$ such that 
$w^{-1}$ has one-line notation $3142$.
We have $\alpha'(w) = (0,1,0,2)$.
The set $\MVPD(w)$ has the following 
three elements:

$$
\begin{tikzpicture}[x=1.5em,y=1.5em,thick,rounded corners, color = blue]
\draw[step=1,gray,ultra thin] (0,0) grid (4,2);
\draw[color=blue] (0,.5)--(2.5, .5)--(2.5, 1.5)--(3.5, 1.5)--(3.5, 2);
\draw[color=blue] (0,1.5)--(1.5, 1.5)--(1.5, 2);
\end{tikzpicture}
\quad\quad\quad
\begin{tikzpicture}[x=1.5em,y=1.5em,thick,rounded corners, color = blue]
\draw[step=1,gray,ultra thin] (0,0) grid (4,2);
\draw[color=blue] (0,.5)--(1.5, .5)--(1.5, 1.5)--(3.5, 1.5)--(3.5, 2);
\draw[color=blue] (0,1.5)--(1.5, 1.5)--(1.5, 2);
\end{tikzpicture}
\quad\quad\quad
\begin{tikzpicture}[x=1.5em,y=1.5em,thick,rounded corners, color = blue]
\draw[step=1,gray,ultra thin] (0,0) grid (4,2);
\draw[color=blue] (0,.5)--(2.5, .5)--(2.5, 1.5)--(3.5, 1.5)--(3.5, 2);
\draw[color=blue] (0,1.5)--(1.5, 1.5)--(1.5, 2);
\node at (2.5,1.5) {$\bullet$};
\end{tikzpicture}
$$

\end{exa}

We describe a bijection from $\PD(w)$
to $\MVPD(w)$.
Take $P \in \PD(w)$ for some $w \in S_n$.
Let $p_1, \cdots, p_k$ be the left-to-right
maxima of $w^{-1}$.
We may remove the pipes $p_1, \cdots, p_k$.
If a $\ptile$ becomes a $\rtile$ after the removal,
we mark it as $\mtile$.
Let $\Phi(P)$ be the resulting tiling.

\begin{exa}
Suppose $n = 8$.
Consider $w \in S_8$ where 
$w^{-1}$ has one-line notation $12547386$.
The left-to-right maxima of $w^{-1}$ are
$1,2,5, 7$ and $8$.
We consider the following $P \in \PD(w)$ 
where pipe $1$, pipe $2$,
pipe $5$, pipe $7$ and pipe $8$
are colored red. Then, the pipes colored blue form $\Phi(P)$ after changing $\rtile$ to $\mtile$ when necessary. 
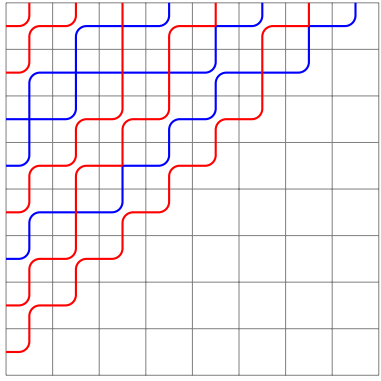
\begin{center}
    \begin{tikzpicture}[x=1.5em,y=1.5em,thick,rounded corners,color = blue]
    \draw[step=1,gray,ultra thin] (0,0) grid (8,8);
     \draw[color=blue, thick] (0,5.5)--(1.5,5.5)--(1.5,7.5)--(3.5,7.5)--(3.5,8);
     \draw[color=blue, thick] (0,4.5)--(0.5,4.5)--(0.5,6.5)--(4.5,6.5)--(4.5,7.5);
     \draw[color=blue, thick] (4.5,7.5)--(5.5,7.5)--(5.5,8);
     \draw[color=blue, thick] (0,2.5)--(0.5,2.5)--(0.5,3.5)--(2.5,3.5)--(2.5,4.5);
     \draw[color=blue, thick](2.5,4.5)--(3.5,4.5)--(3.5,5.5)--(4.5,5.5)--(4.5,6.5)--(6.5,6.5)--(6.5,7.5);
     \draw[color=blue, thick](6.5,7.5)--(7.5,7.5)--(7.5,8);
     \draw[color=red, thick] (0,7.5)--(0.5, 7.5)--(0.5, 8);
     \draw[color=red, thick] (0,6.5)--(0.5, 6.5)--(0.5,7.5)--(1.5, 7.5)--(1.5, 8);
     \draw[color=red, thick] (0,3.5)--(0.5, 3.5)--(0.5,4.5)--(1.5, 4.5)--(1.5, 5.5)--(2.5, 5.5)--(2.5, 8);
     \draw[color=red, thick] (0,1.5)--(0.5, 1.5)--(0.5,2.5)--(1.5, 2.5)--(1.5, 4.5)--(2.5, 4.5);
     \draw[color=red, thick] (2.5, 4.5)--(2.5, 5.5)--(3.5, 5.5)--(3.5, 7.5)--(4.5, 7.5);
     \draw[color=red, thick] (4.5, 7.5)--(4.5, 8);
     \draw[color=red, thick] (0,0.5)--(0.5, 0.5)--(0.5,1.5)--(1.5, 1.5)--(1.5, 2.5)--(2.5, 2.5)--(2.5, 3.5)--(3.5, 3.5)--(3.5, 4.5)--(4.5, 4.5)--(4.5, 5.5)--(5.5, 5.5)--(5.5, 7.5)--(6.5, 7.5);
      \draw[color=red, thick] (6.5, 7.5)--(6.5, 8);
    \end{tikzpicture}
\end{center}
Readers may check that $\Phi(P)$ would be the MVPD in Example~\ref{E: MVPD}. We will later show that $P$ must not have a crossing of two red pipes, which would induce a way to reconstructed $P$ when given only the blue pipes.
\end{exa}

\begin{prop}
\label{P: PD->MVPD}
The map $\Phi$
is a bijection
from $\PD(w)$ to $\MVPD(w)$
that preserves $\wty(\cdot)$
\end{prop}
\begin{proof}
Take any $P\in \PD(w)$ and consider $\Phi(P)$. 
First, if both $P$ and $\Phi(P)$
have pipe $p$,
then it travels the same in $P$ and $\Phi(P)$.
We now check $\Phi(P) \in \MVPD(w)$.
\begin{itemize}
\item We make sure $\Phi(P)$ has no $\vtile$. Suppose to the contrary that $\Phi(P)$ at $(i,j)$ is a $\vtile$, then $P$ must have a $\ptile$ at $(i,j)$. Let pipe $p$ (resp. $q$) be the pipe going horizontally (resp. vertically) in this tile. Then we know pipe $p$ is removed by $\Phi$, so $p$ is a left-to-right maximum of $w^{-1}$. However, since pipe $p$ and pipe $q$ crossed in this tile, we know $p < q$ and $q$ appears on the left of $p$ in the one-line notation of $w^{-1}$, contradicting to $p$ being a left-to-right maximum of $w^{-1}$.
\item Assume pipe $p$ has $\mtile$ at $(i,j)$ of $\Phi(P)$, we check pipe $p$ has a $\htile$ before.
We know $P$ has a $\ptile$ at $(i,j)$. Let pipe $q$ be the other pipe in $(i,j)$ of $P$, so this pipe is removed by $\Phi$. We know pipe $p$ and pipe $q$ already crossed before $(i,j)$ in $P$, say at $(i', j')$.
Then after removing pipe $q$, the $(i', j')$ becomes $\htile$ in $\Phi(P)$. 
\end{itemize}

We have checked $\Phi(P)$ is a valid MVPD. 
Clearly, $\Phi(P)$ has column-to-row code $\alpha'(w)$,
so it is in $\MVPD(w)$. 
We check $\Phi$ preserves $\wty(\cdot)$. 
Take a $\ptile$ in $P$. We check it becomes 
a weighty tile in $\Phi(P)$.
Say pipe $p$ exits from the top and pipe $q$ exits from the right
of this $\ptile$.
By Remark~\ref{R: red pipe no x},
it is impossible that both pipe $p$ and pipe $q$
are removed by $\Phi$,
so this $\ptile$ will not become a $\btile$.
It also cannot become a $\jtile$:
If so, we know $q$ is a left-to-right maximum in $w^{-1}$,
$q < p$, and $q$ ends up on the right
of $p$ in $w^{-1}$. 
This is a contradiction.
It is also obvious that this $\ptile$ cannot be mapped to $\rtile$ or $\bumptile$ by the rules of $\Phi$. Therefore $\Phi$ maps $\ptile$ to weighty tiles. On the other hand, for any $\bumptile$ in $P$, they cannot be mapped to $\htile$, $\ptile$, or $\mtile$ by the rules of $\Phi$. 

Thus, $\Phi$ is a $\wty(\cdot)$ preserving map from $\PD(w)$ to $\MVPD(w)$. It remains to construct its inverse. 
Take $M \in \MVPD(w)$.
We change the cell $(i,j)$ based on the following:
\begin{itemize}
\item If it is $\mtile$ or $\htile$, it becomes $\ptile$.
\item If it is $\jtile$ with $i + j \leq n$, we turn it into $\bumptile$.
If it is $\jtile$ with $i + j = n + 1$, we keep it as $\jtile$.
\item If it is $\rtile$, we turn it into $\bumptile$.
\item Finally, suppose it is $\btile$. 
If $i + j < n$, we turn it into $\bumptile$.
If $i + j = n$, we turn it into $\jtile$.
\end{itemize}

Clearly, we obtain a pipe dream $P$ by adding pipes to each tile. 
We claim for each pipe in $M$,
it goes in the same way in both $P$ and $M$.
In addition, the added pipes in $P$ belong to pipes
which do not exist in $M$. 
We prove by induction on the tiles from bottom to top, 
and left to right in each row.
Consider the tile $(i,j)$
\begin{itemize}
\item If $(i,j)$ is a $\htile$ containing pipe $p$ in $M$,
it becomes a $\ptile$ in $P$.
We need to verify that pipe $p$ goes horizontally in $(i,j)$ of $P$.
Let pipe $q$ be the pipe entering from the bottom
of $(i,j)$ in $P$,
so pipe $q$ does not exist in $P$.
Assume toward contradiction that $p$ does not go horizontally in $(i,j)$.
Then pipe $p$ and pipe $q$ have already crossed,
where the pipe $p$ travels vertically. 
Then the corresponding cell in $M$ would be a $\vtile$,
which is impossible. 
\item If $(i,j)$ is a $\mtile$ containing pipe $p$ in $M$,
it becomes a $\ptile$ in $P$.
We need to verify that pipe $p$ does not go vertically in $(i,j)$ of $P$.
Let pipe $q$ be the pipe entering from the left
of $(i,j)$ in $P$,
so pipe $q$ does not exist in $P$.
We need to show pipe $p$ and $q$ have crossed before. 
Since $(i,j)$ is $\mtile$ in $M$, 
we may find a $\htile$ containing pipe $p$ under row $i$.
In $P$, it becomes a $\ptile$ where pipe $p$
crosses with some added pipe, say pipe $t$.
If $t = q$, we are done. 
Otherwise, we know the added pipes cannot cross. 
Thus, 
the three pipes enter row $i$
with the order $t, q, p$ from left to right. 
By Lemma~\ref{L: abc},
we know pipe $q$ and $p$ have crossed. 
\item The other cases of $(i,j)$ are straightforward to check. 
\end{itemize}

Say $P \in PD(u)$.
The claim above says
the $k\textsuperscript{th}$ entry of $\alpha'(w)$,
if non-zero, agrees with $u^{-1}(k)$. 
Since the added pipes are not crossing,
we know $u^{-1}(k)$ is obtained from $\alpha'(w)$
by turning $0$s into the missing numbers
in increasing order, which yields $w^{-1}$.
Thus, $u = w$ and the map defined above
sends $\MVPD(w)$ to $\PD(w)$.
It is clearly the inverse of $\Phi$.
\end{proof}

\begin{cor}
\label{C: MVPD}
For $w \in S_n$, we have
\begin{align*}
\fG_w(\textbf{x}) & = \sum_{M \in \MVPD(w)}  (-1)^{|\wty(M)| - \ell(w)}\wt_M(\x), \\
\fG_w(\textbf{x}, \textbf{y}) & = \sum_{M \in \MVPD(w)} (-1)^{|\wty(M)| - \ell(w)} \wt_M(\x, \y).  
\end{align*}    
\end{cor}
\begin{proof}
Follows from~(\ref{EQ: Gro formulas}) and Proposition~\ref{P: PD->MVPD}.
\end{proof}

\section{Bumpless Vertical-less Pipe Dreams}
\label{S: BVPD}

\subsection{Describing the BVPD formula}
We introduce \definition{bumpless vertical-less pipe dreams} (BVPD).
They consist of tilings where the following five types of tiles are placed
$$
\btile,\quad \quad \htile, \quad \quad \ptile, \quad \quad \jtile, \quad \quad \rtile
$$
on an $n \times (n-1)$ grid.
The pipes of a BVPD enter from the left edge
and exit from the top edge.
We trace pipes from left to top.
For each $\ptile$,
we trace the pipes in the same way
as PDs and MVPDs.
We name the pipe entering from row $p$
as pipe $p$.
When drawing a BVPD, we omit blank rows
on the bottom and blank
columns to the right. 
\begin{exa}
\label{E: BVPD}
Let $n = 6$.
The following is a BVPD
$$
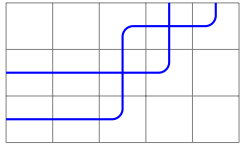
\begin{tikzpicture}[x=1.5em,y=1.5em,thick,rounded corners,color = blue]
\draw[step=1,gray,ultra thin] (0,0) grid (5,3);
\draw[color=blue, thick] (0,0.5)--(2.5, 0.5)--(2.5, 2.5)--(4.5, 2.5)--(4.5, 3);
\draw[color=blue, thick] (0,1.5)--(3.5, 1.5)--(3.5, 3);
\end{tikzpicture}
$$
with pipe $2$ going to column $5$
and pipe $3$ going to column $4$.
\end{exa}

The \definition{column-to-row code}
of a BVPD is a sequence of $n-1$ numbers,
defined similarly as that of a MVPD. 
The column-to-row code of the BVPD in 
Example~\ref{E: BVPD} is $(0,0,0,3,2)$.

Take $w \in S_n$ be inverse fireworks. 
We obtain a sequence $\alpha(w)$ as follows.
We start from the sequence $(w^{-1}(1), \cdots, w^{-1}(n))$
and set the first number in each decreasing run 
to be $0$.
Then $\alpha(w)$ is obtained 
by removing the first entry. 
Notice that $\alpha(w)$
can be obtained from $\alpha'(w)$
by removing the first entry.
Let $\BVPD(w)$ be the set of
all BVPDs with column-to-row code $\alpha(w)$.

\begin{exa}
\label{E: BVPD2}
Say $n = 6$ and $w$ has one-line notation 
$165234$.
Thus, $w^{-1}$ has one-line notation
$145632$ where $1, 4, 5$ and $6$
are the first numbers in the decreasing runs. 
We have $\alpha(w) = (0,0,0,3,2)$,
so $\BVPD(w)$ consists of all
BVPDs with pipe $2$ going to column $5$,
pipe $3$ going to column $4$,
and no other pipes.
There are six such BVPDs:
$$
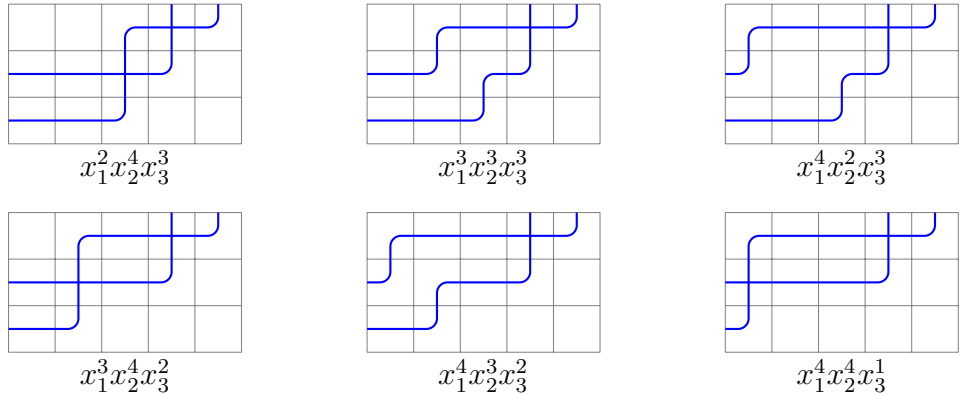
\begin{tikzpicture}[x=1.5em,y=1.5em,thick,rounded corners,color = blue]
\draw[step=1,gray,ultra thin] (0,0) grid (5,3);
\draw[color=blue, thick] (0,0.5)--(2.5, 0.5)--(2.5, 2.5)--(4.5, 2.5)--(4.5, 3);
\draw[color=blue, thick] (0,1.5)--(3.5, 1.5)--(3.5, 3);
\node[color=black] at (2.5,-0.5) {$x_1^2x_2^4x_3^3$};
\end{tikzpicture}
\quad\quad\quad\quad
\begin{tikzpicture}[x=1.5em,y=1.5em,thick,rounded corners,color = blue]
\draw[step=1,gray,ultra thin] (0,0) grid (5,3);
\draw[color=blue, thick] (0,0.5)--(2.5, 0.5)--(2.5, 1.5)--(3.5, 1.5)--(3.5, 3);
\draw[color=blue, thick] (0,1.5)--(1.5, 1.5)--(1.5, 2.5)--(4.5, 2.5)--(4.5, 3);
\node[color=black] at (2.5,-0.5) {$x_1^3x_2^3x_3^3$};
\end{tikzpicture}
\quad\quad\quad\quad
\begin{tikzpicture}[x=1.5em,y=1.5em,thick,rounded corners,color = blue]
\draw[step=1,gray,ultra thin] (0,0) grid (5,3);
\draw[color=blue, thick] (0,0.5)--(2.5, 0.5)--(2.5, 1.5)--(3.5, 1.5)--(3.5, 3);
\draw[color=blue, thick] (0,1.5)--(0.5, 1.5)--(0.5, 2.5)--(4.5, 2.5)--(4.5, 3);
\node[color=black] at (2.5,-0.5) {$x_1^4x_2^2x_3^3$};
\end{tikzpicture}
$$
$$
\begin{tikzpicture}[x=1.5em,y=1.5em,thick,rounded corners,color = blue]
\draw[step=1,gray,ultra thin] (0,0) grid (5,3);
\draw[color=blue, thick] (0,0.5)--(1.5, 0.5)--(1.5, 2.5)--(4.5, 2.5)--(4.5, 3);
\draw[color=blue, thick] (0,1.5)--(3.5, 1.5)--(3.5, 3);
\node[color=black] at (2.5,-0.5) {$x_1^3x_2^4x_3^2$};
\end{tikzpicture}
\quad\quad\quad\quad
\begin{tikzpicture}[x=1.5em,y=1.5em,thick,rounded corners,color = blue]
\draw[step=1,gray,ultra thin] (0,0) grid (5,3);
\draw[color=blue, thick] (0,0.5)--(1.5, 0.5)--(1.5, 1.5)--(3.5, 1.5)--(3.5, 3);
\draw[color=blue, thick] (0,1.5)--(0.5, 1.5)--(0.5, 2.5)--(4.5, 2.5)--(4.5, 3);
\node[color=black] at (2.5,-0.5) {$x_1^4x_2^3x_3^2$};
\end{tikzpicture}
\quad\quad\quad\quad
\begin{tikzpicture}[x=1.5em,y=1.5em,thick,rounded corners,color = blue]
\draw[step=1,gray,ultra thin] (0,0) grid (5,3);
\draw[color=blue, thick] (0,0.5)--(0.5, 0.5)--(0.5, 2.5)--(4.5, 2.5)--(4.5, 3);
\draw[color=blue, thick] (0,1.5)--(3.5, 1.5)--(3.5, 3);
\node[color=black] at (2.5,-0.5) {$x_1^4x_2^4x_3^1$};
\end{tikzpicture}
$$
\end{exa}

Finally, define the \definition{weighty tiles}
of a BVPD $B$, denoted as $\wty(B)$,
as the set of $(i,j)$ that 
is $\ptile$, $\jtile$ or $\htile$ in $B$.
Let the \definition{weight} of $B$, denoted as $\wt_B(\x)$,
be the monomial $\Pi_{(i,j) \in \wty(B)} \: x_i$.
We write the weight of each BVPD under itself in Example~\ref{E: BVPD2}.

\begin{thm}
\label{T: BVPD formula}
For inverse fireworks $w$, we have
$$
\topGro_w(\x) = \sum_{B \in \BVPD(w)} \wt_B(\x).
$$
\end{thm}
Continuing on Example~\ref{E: BVPD2}.
If $w$ has one-line notation $165234$,
then $$\topGro_w(\x) = x_1^2x_2^4x_3^3 + x_1^3x_2^3x_3^3 + x_1^4x_2^2x_3^3 + x_1^3x_2^4x_3^2 + x_1^4x_2^3x_3^2 + x_1^4x_2^4x_3^1.$$

Let $\widehat{\PD(w)}$ be the subset of the $\PD(w)$ with maximal
number of $\ptile$.
\begin{thm}
\label{T: BVPD to PD}
For $w$ inverse fireworks,
there exists a bijection $\Psi$ from 
$\BVPD(w)$ to $\widehat{\PD(w)}$
that 
preserves the positions of weighty tiles.
\end{thm}

Intuitively speaking, for $B \in \BVPD(w)$,
$\Psi(B)$ is the pipe dream 
that has a $\ptile$ at $(i,j)$ when $(i,j-1) \in \wty(B)$ or when $B$ has a pipe entering from row $i$ and $j=1$.
This result also characterizes
the pipe dreams of $w$ with the maximal number of $\ptile$ when $w$ is inverse fireworks. 

\begin{exa}
\label{E: Psi}
Continuing from Example~\ref{E: BVPD2}, the following is $B \in \BVPD(w)$ for $w = 165234$ along with its image under $\Psi$, which we define later.
$$
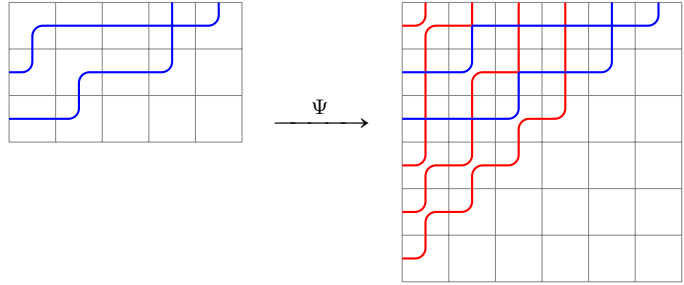
\begin{tikzpicture}[x=1.5em,y=1.5em,thick,rounded corners,color = blue]
\draw[step=1,white,ultra thin] (0,-3) grid (5,0);
\draw[step=1,gray,ultra thin] (0,0) grid (5,3);
\draw[color=blue, thick] (0,0.5)--(1.5, 0.5)--(1.5, 1.5)--(3.5, 1.5)--(3.5, 3);
\draw[color=blue, thick] (0,1.5)--(0.5, 1.5)--(0.5, 2.5)--(4.5, 2.5)--(4.5, 3);
\end{tikzpicture}
\quad\raisebox{2cm}{$\xrightarrow{\quad\Psi\quad}$}\quad 
\begin{tikzpicture}[x=1.5em,y=1.5em,thick,rounded corners,color = blue]
\draw[step=1,gray,ultra thin] (-1,-3) grid (5,3);
\draw[color=red, thick] (-1,2.5)--(-0.5,2.5)--(-0.5,3);
\draw[color=red, thick] (-1,-0.5)--(-0.5,-0.5)--(-0.5,2.5)--(0.5,2.5);
\draw[color=red, thick] (0.5,2.5)--(0.5,3);
\draw[color=red, thick] (-1,-1.5)--(-0.5,-1.5)--(-0.5,-0.5)--(0.5,-0.5)--(0.5,1.5)--(1.5,1.5);
\draw[color=red, thick] (1.5,1.5)--(1.5,3);
\draw[color=red, thick] (-1,-2.5)--(-0.5,-2.5)--(-0.5,-1.5)--(0.5,-1.5)--(0.5,-0.5)--(1.5,-0.5)--(1.5,0.5)--(2.5,0.5)--(2.5,3);
\draw[color=blue, thick] (-1,0.5)--(1.5, 0.5)--(1.5, 1.5);
\draw[color=blue, thick] (1.5, 1.5)--(3.5, 1.5)--(3.5, 3);
\draw[color=blue, thick] (-1,1.5)--(0.5, 1.5)--(0.5, 2.5);
\draw[color=blue, thick] (0.5, 2.5)--(4.5, 2.5)--(4.5, 3);
\end{tikzpicture}
$$
Readers may check that the positions of the weighty tiles are preserved, and that $\Psi(B)$ has maximum crossings in $\PD(w)$.
\end{exa}

\subsection{Proofs of Theorem~\ref{T: BVPD formula} and Theorem~\ref{T: BVPD to PD}}

We start with one simple property on the number of weighty tiles in a MVPD.
For $w \in S_n$, define $r(w) := \sum_{i} i-1$
where $i$ ranges over all numbers such that 
the $i\textsuperscript{th}$ number in $\alpha'(w)$ is non-zero.

\begin{lemma}
\label{L: MVPD weight}
Take $M \in \MVPD(w)$.
Let $k$ be the number of $\bumptile$ and $\rtile$ in $M$.
We have $|\wty(M)| = r(w) - k$.
\end{lemma}
\begin{proof}
We first assign each tile in $M$ that is not $\btile$ or $\jtile$
with a pipe. 
Such a tile must be $\htile$, $\ptile$, $\rtile$, $\mtile$ or $\bumptile$.
We simply assign the pipe that exits from the right
of this tile. 

Take an arbitrary pipe $p$
and suppose it goes to column $c_p$. 
In other words, the $c\textsuperscript{th}$ 
number in $\alpha'(w)$ is $p$.
For each column $i$,
we count the number of cells assigned with pipe $p$
in this column:
\begin{itemize}
\item If the pipe $p$ exits column $i$
and goes to column $i+1$ (i.e. $1 \leq i < c_p$),
there is exactly one tile assigned with
pipe $p$ in column $i$.
\item Otherwise (i.e. $i \geq c_p$), 
there is no tile assigned with pipe $p$.
\end{itemize}

Now there are $c_p-1$ tiles
assigned with the pipe $p$.
Let $k_p$ be the number of 
$\rtile$ and $\bumptile$
assigned with pipe $p$.
The number of weighty tiles
assigned with $p$, denoted as $\wt(p)$, is $(c_p - 1) - k_p$.
We have
\begin{align*}
|\wty(M)|  & = \sum_{\text{pipes }p\text{ in }M}\wt(p) = \sum_{\text{pipes }p\text{ in }M} (c_p - 1) - k_p \\ 
& = \sum_{\text{pipes }p\text{ in }M} (c_p - 1) - \sum_{\text{pipes }p\text{ in }M} k_p
= r(w) - k. \qedhere
\end{align*}
\end{proof}

Let $\widehat{\MVPD(w)}$ be the subset of the $\MVPD(w)$ with the maximal
number of weighty tiles.
Recall that $\raj(w)$ is the degree of 
$\fG_w(\x)$, 
so an element of $\widehat{\MVPD(w)}$
has $\raj(w)$ weighty tiles. 
We can describe $\widehat{\MVPD(w)}$ 
of inverse fireworks $w$ as follows. 

\begin{lemma}
\label{L: No r or bump}
Let $w$ be an inverse fireworks permutation. 
Then $\widehat{\MVPD(w)}$ consists of elements in $\MVPD(w)$
without $\rtile$ and $\bumptile$.
\end{lemma}
\begin{proof}
By Lemma~\ref{L: MVPD weight},
it remains to show $r(w)$ is the 
maximal number of weighty tiles of an element 
in $\MVPD(w)$, which is $\raj(w)$.
By Corollary~\ref{C: raj = raj inv},
$\raj(w) = \raj(w^{-1})$.
Since $w^{-1}$ is fireworks,
by Proposition~\ref{P: raj = maj},
$\raj(w^{-1}) = \maj(w^{-1})$.
It remains to check $r(w) = \maj(w^{-1})$.
 
Recall that $r(w) = \sum_{i \in I} (i-1)$,
where $I = \{i: \textrm{$i\textsuperscript{th}$ entry of $\alpha'(w)$ is not $0$}\}$.
In other words, 
$I$ consists of all $i$ such that 
$w^{-1}(i)$ is not a left-to-right maximum
of $w^{-1}$.
Since $w^{-1}$ is fireworks, 
$I$ consists of $i$
such that $w^{-1}(i)$ is not the 
first number in its decreasing run. 
Then we have
\begin{align*}
\{i-1: i \in I \} & = \{j :\textrm{$w^{-1}(j)$ is not the 
last number in its decreasing run} \}\\
& = \{j : w^{-1}(j) > w^{-1}(j+1)\}.    
\end{align*}

Thus, $r(w) = \sum_{i \in I} (i-1) = \sum_{j : w^{-1}(j) > w^{-1}(j+1)} j = \maj(w^{-1})$.
\end{proof}

\begin{cor}
The first column of any $M \in \widehat{\MVPD(w)}$
only consists of $\btile$ and $\htile$.   
\end{cor}
\begin{proof}
Suppose not. 
Say pipe $p$ has a $\jtile$
in column $1$ of $M$.
We know $w^{-1}(1)$
is a left-to-right maximum in $w^{-1}$,
so the first entry in $\alpha'(w)$ is $0$.
In other words, 
pipe $p$ must exit column $1$ 
and enter column $2$.
Find the cell in column $1$ where
pipe $p$ enters from the bottom and exits from the right.
By $M \in \widehat{\MVPD(w)}$ and Lemma~\ref{L: No r or bump},
this cell can only
be $\ptile$ where the two pipes have already crossed or $\mtile$.
It cannot be a $\ptile$ since pipe $p$ has not crossed
with the pipe entering from the left.
It cannot be a $\mtile$ since pipe $p$ does not have
a $\htile$ before. 
Contradiction.
\end{proof}

Now it remains to establish a bijection
from $\widehat{\MVPD(w)}$ 
to $\BVPD(w)$.
We describe the map $\Phi^{M \rightarrow B}$ as follows. 
Take $M \in \widehat{\MVPD(w)}$, 
we remove its first column and change
all $\mtile$ into $\rtile$, obtaining a tiling $B$.
The inverse of this map, denoted as $\Phi^{B \rightarrow M}$ is also straightforward:
Add a column on the left of $B$ consisting 
of $\htile$ and $\btile$ and
change all $\rtile$ in $B$ into $\mtile$.

\begin{prop}
\label{P: MVPD to BVPD}
The maps $\Phi^{M \rightarrow B}$
and $\Phi^{B \rightarrow M}$ are bijections
between $\widehat{\MVPD(w)}$ and \\
$\BVPD(w)$
that preserve $\wty(\cdot)$.
\end{prop}
\begin{proof}
Say $\Phi^{M \rightarrow B}$ sends 
$M \in \widehat{\MVPD(w)}$ to $B$.
Since $M$ has neither $\vtile$ nor $\bumptile$,
$B$ is a BVPD.
Recall that $\alpha(w)$
is obtained from $\alpha'(w)$
by removing the first $0$.
Since $M$ has column-to-row code $\alpha'(w)$,
we know $B$ has column-to-row code $\alpha(w)$,
so $B \in \BVPD(w)$.
The two maps are clearly inverses of each other.
To show the bijections preserve $\wty(\cdot)$,
we present the following example. 
\end{proof}

\begin{exa}
\label{E: MVPD -> BVPD} The left diagram is $M \in \widehat{\MVPD(w)}$ and the right diagram is $\Phi^{M \rightarrow B}(M) = B \in \BVPD(w)$. 
Their weighty tiles (highlighted yellow)
agree.
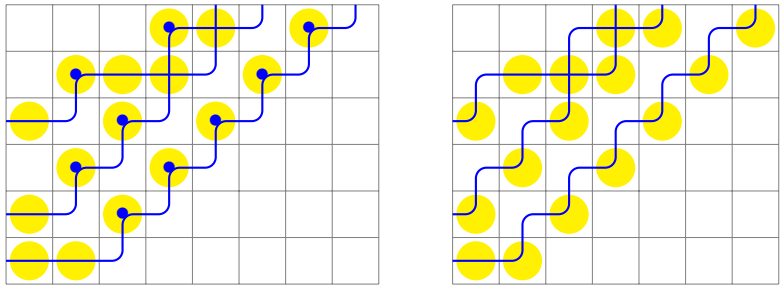
\begin{center}
    \begin{tikzpicture}[x=1.5em,y=1.5em,thick,rounded corners,color = blue]
    \draw[step=1,gray,ultra thin] (0,0) grid (8,6);
    \filldraw [yellow] (0.5,0.5) circle (7pt);
    \filldraw [yellow] (1.5,0.5) circle (7pt);
    \filldraw [yellow] (2.5,1.5) circle (7pt);
    \filldraw [yellow] (3.5,2.5) circle (7pt);
    \filldraw [yellow] (4.5,3.5) circle (7pt);
    \filldraw [yellow] (5.5,4.5) circle (7pt);
    \filldraw [yellow] (6.5,5.5) circle (7pt);
    \filldraw [yellow] (0.5,1.5) circle (7pt);
    \filldraw [yellow] (1.5,2.5) circle (7pt);
    \filldraw [yellow] (2.5,3.5) circle (7pt);
    \filldraw [yellow] (3.5,4.5) circle (7pt);
    \filldraw [yellow] (4.5,5.5) circle (7pt);
    \filldraw [yellow] (0.5,3.5) circle (7pt);
    \filldraw [yellow] (1.5,4.5) circle (7pt);
    \filldraw [yellow] (2.5,4.5) circle (7pt);
    \filldraw [yellow] (3.5,5.5) circle (7pt);
     \draw[color=blue, thick] (0,0.5)--(2.5,0.5)--(2.5,1.5)--(3.5,1.5)--(3.5,2.5)--(4.5,2.5)--(4.5,3.5)--(5.5,3.5)--(5.5,4.5)--(6.5,4.5)--(6.5,5.5)--(7.5,5.5)--(7.5,6);
     \draw[color=blue, thick] (0,1.5)--(1.5,1.5)--(1.5,2.5)--(2.5,2.5)--(2.5,3.5)--(3.5,3.5)--(3.5,5.5)--(5.5,5.5)--(5.5,6);
     \draw[color=blue, thick] (0,3.5)--(1.5,3.5)--(1.5,4.5)--(4.5,4.5)--(4.5,6);
     \node at (6.5,5.5) {$\bullet$};
     \node at (5.5,4.5) {$\bullet$};
     \node at (4.5,3.5) {$\bullet$};
     \node at (3.5,2.5) {$\bullet$};
     \node at (2.5,1.5) {$\bullet$};
     \node at (1.5,2.5) {$\bullet$};
     \node at (2.5,3.5) {$\bullet$};
     \node at (1.5,4.5) {$\bullet$};
     \node at (3.5,5.5) {$\bullet$};
    \end{tikzpicture}
    \quad\quad 
    \begin{tikzpicture}[x=1.5em,y=1.5em,thick,rounded corners,color = blue]
    \filldraw [yellow] (1.5,0.5) circle (7pt);
    \filldraw [yellow] (2.5,0.5) circle (7pt);
    \filldraw [yellow] (3.5,1.5) circle (7pt);
    \filldraw [yellow] (4.5,2.5) circle (7pt);
    \filldraw [yellow] (5.5,3.5) circle (7pt);
    \filldraw [yellow] (6.5,4.5) circle (7pt);
    \filldraw [yellow] (7.5,5.5) circle (7pt);
    \filldraw [yellow] (1.5,1.5) circle (7pt);
    \filldraw [yellow] (2.5,2.5) circle (7pt);
    \filldraw [yellow] (3.5,3.5) circle (7pt);
    \filldraw [yellow] (4.5,4.5) circle (7pt);
    \filldraw [yellow] (5.5,5.5) circle (7pt);
    \filldraw [yellow] (1.5,3.5) circle (7pt);
    \filldraw [yellow] (2.5,4.5) circle (7pt);
    \filldraw [yellow] (3.5,4.5) circle (7pt);
    \filldraw [yellow] (4.5,5.5) circle (7pt);
    \draw[step=1,gray,ultra thin] (1,0) grid (8,6);
     \draw[color=blue, thick] (1,0.5)--(2.5,0.5)--(2.5,1.5)--(3.5,1.5)--(3.5,2.5)--(4.5,2.5)--(4.5,3.5)--(5.5,3.5)--(5.5,4.5)--(6.5,4.5)--(6.5,5.5)--(7.5,5.5)--(7.5,6);
     \draw[color=blue, thick] (1,1.5)--(1.5,1.5)--(1.5,2.5)--(2.5,2.5)--(2.5,3.5)--(3.5,3.5)--(3.5,5.5)--(5.5,5.5)--(5.5,6);
     \draw[color=blue, thick] (1,3.5)--(1.5,3.5)--(1.5,4.5)--(4.5,4.5)--(4.5,6);
    \end{tikzpicture}
\end{center}
\end{exa}

Now we prove the main results of this section.
\begin{proof}[Proof of Theorem~\ref{T: BVPD formula}]
By Corollary~\ref{C: MVPD},
$\topGro_{w}(\x) = 
\sum_{M \in \widehat{\MVPD(w)}} \wt_{M}(\x)$.
Then by Proposition~\ref{P: MVPD to BVPD},
$\sum_{M \in \widehat{\MVPD(w)}} \wt_{M}(\x)
= \sum_{B \in \BVPD(w)} \wt_{B}(\x)$.
\end{proof}

\begin{proof}[Proof of Theorem~\ref{T: BVPD to PD}]
Take $B \in \BVPD(w)$.
To obtain $P \in \widehat{\PD(w)}$,
we simply apply $\Phi^{B \rightarrow M}$ to $B$,
followed by the bijection from $\MVPD(w)$ 
to $\PD(w)$.
Both maps preserve $\wty(\cdot)$,
so $\wty(B) = \wty(P)$.
\end{proof}

\section{Proof of Conjecture~\ref{C: main} for 
inverse fireworks permutations}
\label{S: Proving conjecture}

In this section, we prove
Conjecture~\ref{C: main} for inverse fireworks $w$.
Our approach is constructive: 
For $M \in \MVPD(w) \setminus \widehat{\MVPD(w)}$,
we construct $M'$ such that 
$\wt_M(\x) x_i = \wt_{M'}(\x)$ for some $i$
using ``droop moves''. 
In Section~\ref{SS: Properties of MVPD},
we develop some
general properties of MVPDs 
and define droop moves on MVPDs.
Then we give the construction in Section~\ref{SS: Construction}.

\subsection{Properties of MVPD}
\label{SS: Properties of MVPD}
Let $w \in S_n$ be an arbitrary permutation
in this section.
We start with two observations on $\MVPD(w)$.
\begin{lemma}
\label{L: have horizontal}
Take $M \in \MVPD(w)$.
For every pipe, 
we can find a $\htile$ in $M$
containing that pipe. 
\end{lemma}
\begin{proof}
Consider the pipe from row $r$ of $M$.
We know $r$ appears in $\alpha'(w)$, 
so $r$ is not a left-to-right maximum
in $w^{-1}$.
So there must exist $m>r$, a left-to-right maximum in $w^{-1}$, on the left
of $r$ in $w^{-1}$.
In the pipe dream corresponding to $M$,
there must be a $\ptile$
where the pipe from row $r$ 
goes from left to right
and the pipe from row $m$
goes from bottom to top.
To obtain $M$ from this pipe dream,
we remove the pipe from row $m$,
so this tile becomes a $\htile$.
\end{proof}

We say a $\ptile$ of $M$ is a \definition{real crossing} if its two pipes really cross in it (i.e. the pipe
entering from the bottom exits from top). 
Otherwise, we say the $\ptile$ is a \definition{fake crossing}. 

\begin{lemma}
\label{L: enclosed space}
Take $M \in \MVPD(w)$.
Say pipe $p$ and pipe $q$ have a real crossing
in $(i,j)$ and a fake crossing in $(i', j')$.
We consider the region enclosed by
the two pipes from $(i,j)$ to $(i', j')$.
For any pipe that appears in this region,
it must cross with both pipe $p$ and pipe $q$.
\end{lemma}
\begin{proof}
For a pipe $t$ to enter or exit this region, 
it must cross with pipe $p$ or pipe $q$.
Since two pipes cannot cross more than once,
pipe $t$ must cross both pipe $p$ and pipe $q$.
\end{proof}

Next, we define the \definition{droop moves}
on MVPDs, which look similar
to the droop moves on bumpless pipe dreams introduced
in~\cite{LLS}.
\begin{defn}
\label{D: droop}
Take $M \in \MVPD(w)$.
We define $\droop_{(i,j)}(M)$
if the following are all satisfied
\begin{itemize}
\item The tile $(i,j)$ contains a pipe
that enters from the bottom and exits from 
the right 
(i.e. It is a 
$\rtile$, $\mtile$, $\bumptile$ or a fake crossing).
\item The tile $(i,j + 1)$ is a $\htile$.
\item Let $i' > i$ be the smallest 
such that $(i', j)$ is not $\ptile$.
Then $(i', j)$ is a $\jtile$.
\end{itemize}

For each $i < r < i'$,
we know $(r, j)$ is a $\ptile$.
A simple induction would imply
that $(r, j+1)$ has no pipe entering
from the bottom.
Thus, $(r, j+1)$ is $\htile$
and $(i', j+1)$ is $\rtile$, $\mtile$ or $\btile$.
The operation $\droop(\cdot)$ does the following
to column $j$ and $j+1$ between row $i$ and row $i'$.
\begin{itemize}
\item Change $(i, j)$ from $\bumptile$ or fake crossing 
to $\jtile$. 
Change $(i, j)$ from $\rtile$ or $\mtile$ to $\btile$.
\item Change $(i, j+1)$ from $\htile$ to $\rtile$.
\item For $i < r < i'$,
change $(r, j)$ from $\ptile$ to $\htile$
and change $(r, j+1)$ from $\htile$ to $\ptile$.
\item Change $(i', j)$ from $\jtile$ to $\htile$.
\item Change $(i', j + 1)$ from $\rtile$ or $\mtile$ 
to $\bumptile$.
Change $(i', j + 1)$ from $\btile$ to $\jtile$.
\end{itemize}

We also define $\droop_{(i,j)}'(M)$
on such $(i,j)$ and $M$.
It first performs $\droop_{(i,j)}$.
Then notice that the pipe in $(i, j+1)$
of $\droop_{(i,j)}(M)$ must have 
a $\htile$ in $(r, j)$ for some 
$i < r \leq i'$.
We may mark the pipe in $(i, j+1)$,
obtaining a valid MVPD $\droop_{(i,j)}'(M)$.
\end{defn}

\begin{exa}
We give two examples of the effect of 
$\droop_{(i,j)}'$ and $\droop_{(i,j)}$.
$$
\begin{tikzpicture}[x=1.5em,y=1.5em,thick,rounded corners,color = black]
\draw[step=1,gray,ultra thin] (0,0) grid (2,5);
\draw[color=blue, thick] (0,1.5)--(2,1.5);
\draw[color=blue, thick] (0,2.5)--(2,2.5);
\draw[color=blue, thick] (0,3.5)--(2,3.5);
\draw[color=blue, thick] (0,0.5)--(0.5,0.5)--(0.5,4.5)--(2,4.5);
\node at (-0.5,4.5) {$i$};
\node at (-0.5,0.5) {$i'$};
\node at (0.5,5.5) {$j$};
\end{tikzpicture}
\quad\raisebox{1.45cm}{$\xrightarrow{\droop_{(i,j)}'}$}\quad 
\begin{tikzpicture}[x=1.5em,y=1.5em,thick,rounded corners,color = black]
\draw[step=1,gray,ultra thin] (0,0) grid (2,5);
\draw[color=blue, thick] (0,1.5)--(2,1.5);
\draw[color=blue, thick] (0,2.5)--(2,2.5);
\draw[color=blue, thick] (0,3.5)--(2,3.5);
\draw[color=blue, thick] (0,0.5)--(1.5,0.5)--(1.5,4.5)--(2,4.5);
\node[color=blue] at (1.5,4.5) {$\bullet$};
\node at (-0.5,4.5) {$i$};
\node at (-0.5,0.5) {$i'$};
\node at (0.5,5.5) {$j$};
\end{tikzpicture}
\quad\quad\quad\quad
\begin{tikzpicture}[x=1.5em,y=1.5em,thick,rounded corners,color = black]
\draw[step=1,gray,ultra thin] (0,0) grid (2,5);
\draw[color=blue, thick] (0,1.5)--(2,1.5);
\draw[color=blue, thick] (0,2.5)--(2,2.5);
\draw[color=blue, thick] (0,3.5)--(2,3.5);
\draw[color=blue, thick] (0,0.5)--(0.5,0.5)--(0.5,4.5)--(2,4.5);
\draw[color=blue, thick] (0,4.5)--(.5, 4.5)--(.5, 5);
\draw[color=blue, thick] (1.5, 0)--(1.5, 0.5)--(2,.5);
\node[color=blue] at (1.5,0.5) {$\bullet$};
\node at (-0.5,4.5) {$i$};
\node at (-0.5,0.5) {$i'$};
\node at (0.5,5.5) {$j$};
\end{tikzpicture}
\quad\raisebox{1.45cm}{$\xrightarrow{\droop_{(i,j)}}$}\quad 
\begin{tikzpicture}[x=1.5em,y=1.5em,thick,rounded corners,color = black]
\draw[step=1,gray,ultra thin] (0,0) grid (2,5);
\draw[color=blue, thick] (0,1.5)--(2,1.5);
\draw[color=blue, thick] (0,2.5)--(2,2.5);
\draw[color=blue, thick] (0,3.5)--(2,3.5);
\draw[color=blue, thick] (0,0.5)--(1.5,0.5)--(1.5,4.5)--(2,4.5);
\draw[color=blue, thick] (0,4.5)--(.5, 4.5)--(.5, 5);
\draw[color=blue, thick] (1.5, 0)--(1.5, 0.5)--(2,.5);
\node at (-0.5,4.5) {$i$};
\node at (-0.5,0.5) {$i'$};
\node at (0.5,5.5) {$j$};
\end{tikzpicture}
$$
\end{exa}

\begin{lemma}
\label{L: Droop}
Take $M \in \MVPD(w)$.
Then $\droop_{(i,j)}(M)$ and $\droop_{(i,j)}'(M)$
are both in $\MVPD(w)$ if they are defined. 
\end{lemma}
\begin{proof}
Let $i' > i$ be the smallest such that $(i',j)$
is not $\ptile$.
We just need to show that the same pipe exits 
from the right edge of $(r, j+1)$
for $i \leq r < i'$ in $M$ and $\droop_{(i,j)}(M)$.
To prove this, we claim: 
For $i < r \leq i'$, 
if a pipe exits
from the top of $(r, j)$ in $M$,
then the same pipe exits from the top 
of $(r, j+1)$ in $\droop_{(i,j)}(M)$.
We prove by induction on $r$.
The base case when $r = i'$ is immediate.
Now suppose $i < r < i'$.
Say pipe $p$ enters $(r,j)$ from the left
and pipe $q$ enters $(r,j)$ from the bottom
in $M$.
Then pipe $\max(p,q)$ exits from the the top
of $(r,j)$. 
The other pipe exits from the right of $(r, j+1)$.
By our inductive hypothesis, 
pipe $p$ enters $(r,j+1)$ from the left
and pipe $q$ enters $(r,j+1)$ from the bottom
in $\droop_{(i,j)}(M)$.
Pipe $\max(p,q)$ exits from the the top
of $(r,j+1)$,
and the other pipe exits from the right.
Our inductive step is finished. 
\end{proof}

Finally, we study a special family of MVPDs.
\begin{defn}
A $M \in \MVPD(w)$ is called
\definition{saturated} if it satisfies both of the following. 
\begin{itemize}
\item For any $\rtile$ in $M$, 
the pipe in it does not have $\htile$ before.
\item For any $\bumptile$,
the two pipes in it do not cross in $M$. 
\end{itemize}

In other words, an $M \in \MVPD(w)$ is not saturated
if we can turn one of its $\bumptile$ to
$\ptile$ or $\rtile$ to $\mtile$
and still remain in $\MVPD(w)$.
\end{defn}

\begin{lemma}
\label{L: no bad + real crossing}
Take a saturated $M \in \MVPD(w)$.
Say a pipe $p$ enters the tile $(i,j)$
from the bottom and exits from the right.
Then $(i, j + 1)$ cannot be a real crossing. 
\end{lemma}
\begin{proof}
Suppose there exists such $(i,j)$.
We pick one such $(i,j)$ 
where $i$ is maximal.
Say pipe $p$ enters from the bottom of $(i,j)$
and say it crosses with pipe $q$ in $(i, j+1)$.
We know these two pipes have not crossed
before $(i, j + 1)$.
Moreover, since $M$ is saturated,
there is no $\bumptile$ in $M$
involving pipe $p$ and pipe $q$.
As a conclusion, 
under row $i$,
there is no tile containing both pipe $p$
and pipe $q$.

Find $i' > i$ such that 
pipe $p$ enters on the left edge of $(i', j)$.
We know pipe $p$ goes from bottom to top
of $(r, j)$ for $i < r < i'$.
Say pipe $q$ enters on the left edge of $(i'', j+1)$.
We have $i'' > i'$ since otherwise,
$(i'', j)$ would be a tile containing both pipe $p$
and pipe $q$.
Pipe $q$ goes from bottom to top
of $(i', j+1)$, so $(i', j+1)$ is a real crossing.
Consider the tile $(i', j)$.
It contains pipe $p$ which enters on the left and 
exits on the top. 
Thus, it also must contain a pipe 
that enters from the bottom
and exits on the right.
We reach a contradiction since
we picked the maximal $i$.
\end{proof}

Here is an illustration of the proof 
of Lemma~\ref{L: no bad + real crossing}.
We make pipe $p$ green and pipe $q$ red.
$$
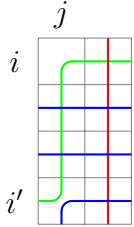
\begin{tikzpicture}[x=1.5em,y=1.5em,thick,rounded corners,color = black]
\draw[step=1,gray,ultra thin] (0,0) grid (2,4);
\draw[color=green, thick] (0,0.5)--(0.5,0.5)--(0.5,3.5)--(2,3.5);
\draw[color=red, thick] (1.5,0)--(1.5,4);
\draw[color=blue, thick] (0.5,0)--(0.5,0.5)--(2,0.5);
\draw[color=blue, thick] (0,1.5)--(2,1.5);
\draw[color=blue, thick] (0,2.5)--(2,2.5);
\node at (-0.5,3.5) {$i$};
\node at (-0.5,0.5) {$i'$};
\node at (0.5,4.5) {$j$};
\end{tikzpicture}
$$

\subsection{Construction}
\label{SS: Construction}
Fix inverse fireworks $w$ throughout this section. 
For each non-maximal marked vertical-less pipe dream $M \in \MVPD(w) \setminus \widehat{\MVPD(w)}$,
our goal is to construct $M'$
such that $\wt_{M'}(\x) = \wt_M(\x) x_i$
for some $i$.
If $M$
is not saturated, we can find the $M'$ easily:
Say $M$ has an $\rtile$ 
and the pipe in it has a $\htile$ before, 
we simply mark the $\rtile$ and obtain $M'$.
Otherwise, say $M$ has a $\bumptile$
where the two pipes in it cross somewhere else. 
We may turn this $\bumptile$ into $\ptile$
and the resulting MVPD is still in $\MVPD(w)$.
It remains to consider saturated $M$. 
Our construction relies 
on the operator $\droop_{(i,j)}'(\cdot)$,
which requires us to find an occurrence of
$\bumphtile$ or $\rhtile$ in $M$.
That is, a $\bumptile$ or $\rtile$
with a $\htile$ immediately on its right.

\begin{lemma}
\label{L: bad horizontal}
Take a saturated $M \in \MVPD(w) \setminus \widehat{\MVPD(w)}$.
In $M$, there exists $\bumphtile$
or $\rhtile$.
\end{lemma}
\begin{proof}
Since $M \notin \widehat{\MVPD(w)}$,
by Lemma~\ref{L: No r or bump},
$M$ must have a $\bumptile$ or $\rtile$.
Let $(i,j)$ be the highest, or one of the highest, such tile.
We prove $(i, j + 1)$ must be $\htile$ by contradiction.
Suppose $(i, j+1)$ is not a $\htile$.
Let pipe $p$ be the pipe
that enters $(i,j)$ from the bottom and exits on the right. 
By Lemma~\ref{L: no bad + real crossing},
$(i, j+1)$ cannot be a real crossing. 
Thus, $(i, j+1)$ can be a fake crossing, 
a $\bumptile$, or a $\jtile$.
In any case, 
pipe $p$ must exit on the top
of $(i, j+1)$.
Then we present two different arguments
based on whether $(i,j)$
is $\bumptile$ or $\rtile$.
Both arguments follow the following three steps: 
\begin{itemize}
\item Step 1: Show pipe $p$ must exit column $j+1$ from the right.
Say it exits from the right edge of $(i', j+1)$ for some $i' < i$.
\item Step 2: We know $(i',j+1)$ cannot be $\rtile$ or $\bumptile$
by how we picked $(i,j)$.
We show $(i', j+1)$ cannot be $\mtile$, so it must be a fake crossing. 
\item Step 3: Find a contradiction. 
\end{itemize}

We start with the case where $(i, j)$ is $\bumptile$. 
Let pipe $q$ be the pipe exiting from the top of $(i,j)$.
Since $M$ is saturated, we know 
pipe $p$ and pipe $q$ never cross in $M$,
so $q < p$.
Now we perform the three steps and 
eventually show pipe $p$ and pipe $q$ must cross,
which would be a contradiction. 
\begin{itemize}
\item Step $1$:
Suppose pipe $p$ does not exit column $j+1$ from the right, then it must exit column $j+1$ from the top. 
Since pipe $p$ and $q$ cannot cross, 
pipe $q$ cannot exit column $j$ from the right.
In other words, $\alpha'(j+1) = p$ and $\alpha'(j) = q$.
Then $w^{-1}(j+1) = p$ and $w^{-1}(j) = q$.
Since $q < p$,
$p$ is actually the first number in its decreasing run
in $w^{-1}$, so $p$ cannot appear in $\alpha(w)$.
We reach a contradiction, so pipe $p$ must exit column $j+1$ from the right.

\item Step $2$:
We know pipe $p$ goes from the bottom to top in $(r, j+1)$
for $i' < r < i$.
Since pipe $q$ cannot cross with pipe $p$,
it must also go from bottom to top in $(r, j)$
for $i' < r < i$.
Thus, pipe $q$ enters $(i', j)$ from the bottom.
The tile $(i', j)$ must have some pipe
exiting from the right. 
Thus, $(i', j + 1)$ has a pipe entering from the left,
so it cannot be $\mtile$.
It must be a fake crossing. 

\item Step $3$:
Let pipe $t$ be the pipe that enters $(i', j+1)$
from the left. 
Since $(i', j+1)$ is a fake crossing, 
pipe $t$ and pipe $p$ must have a real crossing under row $i$.
Then pipe $t$ must exit row $i$ on the left of pipe $p$.
Since pipe $p$ exits row $i$ on column $j+1$
and pipe $q$ exits row $i$ on column $j$,
we know pipe $t$ exits row $i$ on the left 
of column $j$.
Now consider the region enclosed by pipe $p$
and $t$ from their real crossing to $(i', j+1)$.
Pipe $q$ appears in this region.
By Lemma~\ref{L: enclosed space},
pipe $q$ crosses with pipe $p$. 
Contradiction. 
\end{itemize}

Now assume $(i,j)$ is $\rtile$.
By $M$ is saturated,
we know pipe $p$ does not have a $\htile$
before $(i, j)$.
We perform the three steps. 
\begin{itemize}
\item Step $1$:
If pipe $p$ does not exit column $j+1$ from the right,
then it does not have a $\htile$ in $M$,
contradicting Lemma~\ref{L: have horizontal}.
\item Step $2$:
Pipe $p$ still
does not have a $\htile$ before $(i',j+1)$,
so $(i',j+1)$ cannot be $\mtile$.
It must be a fake crossing.
\item Step $3$:
Say $(i', j+1)$ is a fake crossing
between pipe $p$ and pipe $t$. 
Pipe $p$ must have a real crossing under row $i$
where pipe $p$ goes horizontally.
In other words, we can find a real crossing $(r_+, c_+)$ where 
pipe $p$ goes horizontally with $r_+ > i$.
Take the $(r_+, c_+)$ where $c_+$ is maximal.
Thus, from $(r_+, c_+)$ to $(i, j+1)$,
the pipe $p$ is not allowed to travel
horizontally in any tile. 
In other words, if pipe $p$ enters a tile from 
the left, it must exit from the top. 

Next, we argue for $i \leq r \leq r_+$, 
when pipe $p$ exits row $r$,
there is a pipe exiting from the cell on its left, 
which has already crossed with pipe $p$.

We prove our claim by induction.
The base case is when $r = r_+$.
We know $(r_+, c_+)$ is a real crossing.
Since pipe $p$ enters $(r_+, c_+ +1)$ from the left,
it exits row $r_+$ from $(r_+, c_+ +1)$.
Indeed, $(r_+, c_+)$ has a pipe exiting from the top,
which just crossed with pipe $p$.
Now take $i \leq r < r_+$.
Say pipe $p$ enters from the bottom of $(r,c)$.
By our inductive hypothesis, 
another pipe enters $(r,c-1)$ from the bottom, which we label as pipe $s$.
If pipe $p$ goes vertically in $(r,c)$,
pipe $s$ must go vertically in $(r,c-1)$
since if it exits on the right, 
$(r,c)$ would be a fake crossing. 
Now suppose pipe $p$ exits $(r,c)$
on the right. 
Since $(r,c-1)$ has a pipe entering from the bottom,
it must have a pipe exiting from the right.
Then $(r,c)$ can be a fake crossing or a $\bumptile$.
Consider the region enclosed by pipe $t$
and pipe $p$ from their real crossing to $(i', j+1)$.
The other pipe in $(r,c)$ is either pipe $t$,
or lies in this region. 
In either case, it must cross with pipe $p$.
Since $M$ is saturated,
$(r,c)$ is not a $\bumptile$, so it is a fake crossing.
Pipe $p$ exits from the top of $(r,c+1)$
and some pipe that has crossed with it exits from the 
top of $(r,c)$.

Finally, our claim implies when pipe $p$ exits $(i,j+1)$ from the top, there must be a pipe that exits $(i,j)$ from the top.
This contradicts to our assumption that $(i,j)$ is $\rtile$. \qedhere
\end{itemize}
\end{proof}

Now we describe our algorithm.
Take a saturated $M \in \MVPD(w) \setminus \widehat{\MVPD(w)}$.
Lemma~\ref{L: bad horizontal} implies that $M$ must have a 
$\bumphtile$ or $\rhtile$.
We let $(i,j)$ and $(i, j+1)$ be the
lowest such occurrence
where we first maximize $i$,
and then $j$.
We check $\droop_{(i,j)}'(M)$
is defined.
The first two conditions in Definition~\ref{D: droop} are immediate. 
For the last condition,
we let $i' > i$ be the smallest
such that $(i', j)$ is not $\ptile$.
It can be $\bumptile$ or $\jtile$.
Assume it is a $\bumptile$ toward
contradiction. 
Since $(i'-1, j+1)$ is $\htile$,
we know $(i', j+1)$ is also a $\htile$.
This contradicts the maximality of $i$.
Thus, $(i', j)$ is a $\jtile$
and $\droop_{(i,j)} (M)$ 
is well-defined. 

Next, the algorithm computes $\droop_{(i,j)}'(M)$,
which is in $\MVPD(w)$ by Lemma~\ref{L: Droop}.
We compare the weighty tiles of $M$
and $\droop_{(i,j)}'(M)$:
\begin{itemize}
\item The tile $(i, j)$ is not 
weighty in $M$ and $\droop_{(i,j)}'(M)$.
The tile $(i, j+1)$ is weighty
in $M$ and $\droop_{(i,j)}'(M)$.
\item For $i < r < i'$,
the tile $(r, j)$ and $(r, j+1)$ 
are weighty in $M$ and $\droop_{(i,j)}'(M)$.
\item The tile $(i', j)$ is not
weighty in $M$ but 
becomes weighty in $\droop_{(i,j)}'(M)$.
The tile $(i', j+1)$ could be either
weighty or not 
in $M$, but is not weighty in $\droop_{(i,j)}'(M)$.
\end{itemize}

If $\droop_{(i,j)}'(M)$ has one more
weighty tile than $M$, 
we let $M' = \droop_{(i,j)}'(M)$
and terminate.
Then $\wt_{M'}(\x) = \wt_{M}(\x) x_{i'}$.
Otherwise, 
$$
\wty(\droop_{(i,j)}'(M))
= (\wty(M) \setminus \{(i', j+1\}) \cup \{(i', j)\},
$$ 
so $\droop_{(i,j)}'(M)$ and $M$ have the same
number of weighty tiles in each row. 
If $\droop_{(i,j)}'(M)$ is not saturated, 
then we change a $\rtile$ or a $\bumptile$
into a weighty tile and obtain $M'$.
Otherwise, we update the variable $M$
into $\droop_{(i,j)}'(M)$ and repeat the algorithm. 

It remains to show the algorithm eventually terminates.
Let $M_1, M_2, \cdots$ be the MVPDs
in the start of each iteration. 
We know $\wty(M_k)$
is obtained from $\wty(M_{k-1})$
by turning an $(r,c)$ into $(r, c-1)$. 
Thus, the algorithm must terminate.

\begin{exa}
The following is an example of the algorithm.
We start with a saturated $M \in \MVPD(w)\setminus \widehat{\MVPD(w)}$
where $w^{-1}$ has one-line notation $14253$.
We first apply $\droop_{(1,1)}'$ and obtain
another $M_2 \in \MVPD(w)$.
Notice that $M_2$ is also saturated 
and $\wt_M(\x) = \wt_{M_2}(\x)$.
We then apply $\droop_{(2,2)}'$ and obtain $M'$.
Notice that $\wt_{M'}(\x) = \wt_M(\x) x_{3}$.

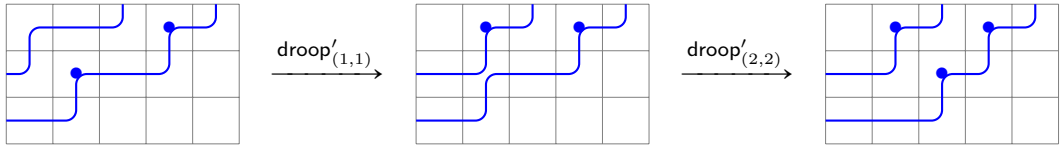
$$
    \begin{tikzpicture}[x=1.5em,y=1.5em,thick,rounded corners,color = blue]
    \draw[step=1,gray,ultra thin] (0,0) grid (5,3);
     \draw[color=blue, thick] (0,0.5)--(1.5,0.5)--(1.5,1.5)--(3.5,1.5)--(3.5,2.5)--(4.5,2.5)--(4.5,3);
     \draw[color=blue, thick] (0,1.5)--(0.5,1.5)--(0.5,2.5)--(2.5,2.5)--(2.5,3);
     \node at (1.5,1.5) {$\bullet$};
     \node at (3.5,2.5) {$\bullet$};
    \end{tikzpicture}
    \quad\raisebox{0.82cm}{$\xrightarrow{\droop_{(1,1)}'}$}\quad 
    \begin{tikzpicture}[x=1.5em,y=1.5em,thick,rounded corners,color = blue]
    \draw[step=1,gray,ultra thin] (0,0) grid (5,3);
    \draw[color=blue, thick] (0,0.5)--(1.5,0.5)--(1.5,1.5)--(3.5,1.5)--(3.5,2.5)--(4.5,2.5)--(4.5,3);
    \draw[color=blue, thick] (0,1.5)--(1.5,1.5)--(1.5,2.5)--(2.5,2.5)--(2.5,3);
    \node at (1.5,2.5) {$\bullet$};
    \node at (3.5,2.5) {$\bullet$};
    \end{tikzpicture}
    \quad\raisebox{0.82cm}{$\xrightarrow{\droop_{(2,2)}'}$}\quad 
    \begin{tikzpicture}[x=1.5em,y=1.5em,thick,rounded corners,color = blue]
    \draw[step=1,gray,ultra thin] (0,0) grid (5,3);
    \draw[color=blue, thick] (0,0.5)--(2.5,0.5)--(2.5,1.5)--(3.5,1.5)--(3.5,2.5)--(4.5,2.5)--(4.5,3);
    \draw[color=blue, thick] (0,1.5)--(1.5,1.5)--(1.5,2.5)--(2.5,2.5)--(2.5,3);
    \node at (1.5,2.5) {$\bullet$};
    \node at (3.5,2.5) {$\bullet$};
    \node at (2.5,1.5) {$\bullet$};
    \end{tikzpicture}
$$
\end{exa}

\section{Acknowledgements}
We thank Brendon Rhoades for reading an earlier version of this paper and giving many useful comments. We thank the anonymous referees for their helpful suggestions to improve this paper.

\bibliographystyle{alpha}
\bibliography{main.bbl}{}

\newcommand{\etalchar}[1]{$^{#1}$}
\begin{thebibliography}{HMMSD22}

\bibitem[BB93]{BB}
Nantel Bergeron and Sara Billey.
\newblock {R}{C}-graphs and {S}chubert polynomials.
\newblock {\em Experimental Mathematics}, 2(4):257--269, 1993.

\bibitem[BJS93]{BJS}
Sara~C Billey, William Jockusch, and Richard~P Stanley.
\newblock Some combinatorial properties of {S}chubert polynomials.
\newblock {\em Journal of Algebraic Combinatorics}, 2(4):345--374, 1993.

\bibitem[CFY24]{CFY}
Yiming Chen, Neil~J.Y. Fan, and Zelin Ye.
\newblock Zero-one {G}rothendieck polynomials.
\newblock {\em arXiv preprint arXiv:2405.05483}, 2024.

\bibitem[CY24]{CY}
Chen-An Chou and Tianyi Yu.
\newblock Constructing maximal pipedreams of double {G}rothendieck polynomials.
\newblock {\em The Electronic Journal of Combinatorics}, 31(3):P3.15, 2024.

\bibitem[DMSD24]{DMS}
Matt Dreyer, Karola Me{\'s}z{\'a}ros, and Avery St~Dizier.
\newblock On the degree of {G}rothendieck polynomials.
\newblock {\em Algebraic Combinatorics}, 7(3):627--658, 2024.

\bibitem[FK94]{FK}
Sergey Fomin and Anatol~N Kirillov.
\newblock Grothendieck polynomials and the {Y}ang-{B}axter equation.
\newblock In {\em Proc. formal power series and alg. comb}, pages 183--190,
  1994.

\bibitem[Haf22]{Haf}
Elena~S Hafner.
\newblock Vexillary {G}rothendieck polynomials via bumpless pipe dreams.
\newblock {\em arXiv preprint arXiv:2201.12432}, 2022.

\bibitem[HMMSD22]{HMMS}
June Huh, Jacob Matherne, Karola M{\'e}sz{\'a}ros, and Avery St.~Dizier.
\newblock Logarithmic concavity of {S}chur and related polynomials.
\newblock {\em Transactions of the American Mathematical Society},
  375(6):4411--4427, 2022.

\bibitem[HMSSD23]{HMSS}
Elena~S Hafner, Karola M{\'e}sz{\'a}ros, Linus Setiabrata, and Avery
  St.~Dizier.
\newblock {M}-convexity of {G}rothendieck polynomials via bubbling.
\newblock {\em arXiv preprint arXiv:2306.08597}, 2023.

\bibitem[KM05]{KM}
Allen Knutson and Ezra Miller.
\newblock Gr{\"o}bner geometry of {S}chubert polynomials.
\newblock {\em Annals of Mathematics}, pages 1245--1318, 2005.

\bibitem[LLS21]{LLS}
Thomas Lam, Seung~Jin Lee, and Mark Shimozono.
\newblock Back stable {S}chubert calculus.
\newblock {\em Compositio Mathematica}, 157(5):883--962, 2021.

\bibitem[LLS23]{LLS2}
Thomas Lam, Seung~Jin Lee, and Mark Shimozono.
\newblock Back stable {K}-theory {S}chubert calculus.
\newblock {\em International Mathematics Research Notices},
  2023(24):21381--21466, 2023.

\bibitem[LS82]{LS:Groth}
Alain Lascoux and Marcel-Paul Sch\"{u}tzenberger.
\newblock Structure de {H}opf de l'anneau de cohomologie et de l'anneau de
  {G}rothendieck d'une vari\'{e}t\'{e} de drapeaux.
\newblock {\em C. R. Acad. Sci. Paris S\'{e}r. I Math.}, 295(11):629--633,
  1982.

\bibitem[MSSD24]{MSS}
Karola M{\'e}sz{\'a}ros, Linus Setiabrata, and Avery St.~Dizier.
\newblock On the support of {G}rothendieck polynomials.
\newblock {\em Annals of Combinatorics}, pages 1--22, 2024.

\bibitem[PSW24]{PSW}
Oliver Pechenik, David~E Speyer, and Anna Weigandt.
\newblock Castelnuovo--mumford regularity of matrix schubert varieties.
\newblock {\em Selecta Mathematica}, 30(4):66, 2024.

\bibitem[RRR{\etalchar{+}}21]{RRRSW}
Jenna Rajchgot, Yi~Ren, Colleen Robichaux, Avery St.~Dizier, and Anna Weigandt.
\newblock Degrees of symmetric {G}rothendieck polynomials and
  {C}astelnuovo-{M}umford regularity.
\newblock {\em Proceedings of the American Mathematical Society},
  149(4):1405--1416, 2021.

\end{thebibliography}
\end{document}